\documentclass[12pt]{article}
\usepackage{amsmath}
\usepackage{amsfonts}
\usepackage{amssymb}

\def\thebibliograph#1#2{\section*{{\normalsize \bf #2}}\list
   {[\arabic{enumi}]}{\settowidth\labelwidth{[#1]}\leftmargin\labelwidth
     \advance\leftmargin\labelsep
     \usecounter{enumi}}
     \def\newblock{\hskip .11em plus .33em minus -.07em}
     \sloppy
     \sfcode`\.=1000\relax}

\newtheorem{theorem}{Theorem}

\newtheorem{definition}{Definition}
\newtheorem{corollary}{Corollary}
\newtheorem{lemma}{Lemma}
\newtheorem{remark}{Remark}
\input{tcilatex}
\begin{document}

\title{Jawerth-Franke embeddings of Herz-type Besov and Triebel-Lizorkin
spaces}
\author{ Douadi Drihem \ \thanks{%
M'sila University, Department of Mathematics, Laboratory of Functional
Analysis and Geometry of Spaces , P.O. Box 166, M'sila 28000, Algeria,
e-mail: \texttt{\ douadidr@yahoo.fr}}}
\date{\today }
\maketitle

\begin{abstract}
In this paper we prove the Jawerth-Franke embeddings of Herz-type Besov and
Triebel-Lizorkin spaces. Moreover, we obtain the Jawerth-Franke embeddings
of Besov and Triebel-Lizorkin spaces equipped with power weights. An
application\ we present new embeddings between Besov and Herz spaces.
\medskip


\end{abstract}

\section{Introduction}


The Herz-type Besov-Triebel-Lizorkin spaces initially appeared in the papers
of J. Xu and D. Yang \cite{XuYang05} and \cite{Xu03}. Several basic
properties were established, such as the Fourier analytical characterisation
and\ lifting properties. When $\alpha =0$ and$\ p=q$ they coincide with the
usual function spaces $F_{p,q}^{s}$.\vskip5pt

The interest in these spaces comes not only from theoretical reasons but
also from their applications to several classical problems in analysis. In 
\cite{LuYang97}, Lu and Yang introduced the Herz-type Sobolev and Bessel
potential spaces. They gave some applications to partial differential
equations. Also in \cite{T11}, Y. Tsutsui, studied the Cauchy problem for
Navier-Stokes equations on Herz spaces\ and weak Herz spaces.

\vskip5pt

The main aim of this paper is to prove the Jawerth-Franke embeddings in $%
\dot{K}_{q}^{\alpha ,p}F_{\beta }^{s}$ and $\dot{K}_{q}^{\alpha ,p}B_{\beta
}^{s}$ spaces, where we use the so-called $\varphi $-transform
characterization in the sense of Frazier and Jawerth. As a consequence, we
present the Jawerth-Franke embeddings of Besov and Triebel-Lizorkin spaces
equipped with power weights. Also, we present new embeddings between Besov
and Herz spaces. All these results generalize the existing classical results
on Besov and Triebel-Lizorkin spaces.

For any $u>0,k\in \mathbb{Z}$ we set $C\left( u\right) =\{x\in \mathbb{R}%
^{n}:u/2\leq \left\vert x\right\vert <u\}$ and $C_{k}=C(2^{k})$. For $x\in 
\mathbb{R}^{n}$ and $r>0$ we denote by $B(x,r)$ the open ball in $\mathbb{R}%
^{n}$ with center $x$ and radius $r$. Let $\chi _{k}$, for $k\in \mathbb{Z}$%
, denote the characteristic function of the set $C_{k}$. The expression $%
f\approx g$ means that $C$ $g\leq f\leq c\,g$ for some independent constants 
$c,C$ and non-negative functions $f$ and $g$. \vskip5pt

\noindent We denote by $\left\vert \Omega \right\vert $ the $n$-dimensional
Lebesgue measure of $\Omega \subseteq \mathbb{R}^{n}$. For any measurable
subset $\Omega \subseteq \mathbb{R}^{n}$ the Lebesgue space $L^{p}(\Omega )$%
, $0<p\leq \infty $ consists of all measurable functions for which $\big\|%
f\mid L^{p}(\Omega )\big\|=\big(\int_{\Omega }\left\vert f(x)\right\vert
^{p}dx\big)^{1/p}<\infty $, $0<p<\infty $ and $\big\|f\mid L^{\infty
}(\Omega )\big\|=\underset{x\in \Omega }{\text{ess-sup}}\left\vert
f(x)\right\vert <\infty $. If $\Omega =\mathbb{R}^{n}$ we put $L^{p}(\mathbb{%
R}^{n})=L^{p}$ and $\big\|f\mid L^{p}(\mathbb{R}^{n})\big\|=\big\|f\big\|%
_{p} $. The symbol $\mathcal{S}(\mathbb{R}^{n})$ is used in place of the set
of all Schwartz functions and we denote by $\mathcal{S}^{\prime }(\mathbb{R}%
^{n})$ the dual space of all tempered distributions on $\mathbb{R}^{n}$. We
define the Fourier transform of a function $f\in \mathcal{S}(\mathbb{R}^{n})$
by $\mathcal{F(}f)(\xi )=\left( 2\pi \right) ^{-n/2}\int_{\mathbb{R}%
^{n}}e^{-ix\cdot \xi }f(x)dx$. Its inverse is denoted by $\mathcal{F}^{-1}f$%
. Both $\mathcal{F}$ and $\mathcal{F}^{-1}$ are extended to the dual
Schwartz space $\mathcal{S}^{\prime }(\mathbb{R}^{n})$ in the usual way.%
\vskip5pt

Let $\mathbb{Z}^{n}$ be the lattice of all points in $\mathbb{R}^{n}$ with
integer-valued components. If $v\in \mathbb{N}_{0}$ and $m=(m_{1},...,m_{n})%
\in \mathbb{Z}^{n}$ we denote $Q_{v,m}$ the dyadic cube in $\mathbb{R}^{n};$%
\begin{equation*}
Q_{v,m}=\{(x_{1},...,x_{n}):m_{i}\leq 2^{v}x_{i}<m_{i}+1,i=1,2,...,n\}.
\end{equation*}%
By $\chi _{v,m}$ we denote the characteristic function of the cube $Q_{v,m}$.

Given two quasi-Banach spaces $X$ and $Y$, we write $X\hookrightarrow Y$ if $%
X\subset Y$ and the natural embedding of $X$ in $Y$ is continuous. We use $c$
as a generic positive constant, i.e.\ a constant whose value may change from
appearance to appearance.\vskip5pt

\section{Function spaces}

We start by recalling the definition and some of the properties of the
homogenous Herz spaces $\dot{K}_{q}^{\alpha ,p}$.\vskip5pt

\begin{definition}
\label{Herz-spaces} \textit{Let }$\alpha \in \mathbb{R},0<p,q\leq \infty $%
\textit{. The homogeneous Herz space }$\dot{K}_{q}^{\alpha ,p}$\textit{\ is
defined by }%
\begin{equation*}
\dot{K}_{q}^{\alpha ,p}=\{f\in L_{\mathrm{loc}}^{q}(\mathbb{R}^{n}\setminus
\{0\}):\big\|f\big\|_{\dot{K}_{q}^{\alpha ,p}}<\infty \},
\end{equation*}%
\textit{where }%
\begin{equation*}
\big\|f\big\|_{\dot{K}_{q}^{\alpha ,p}}=\Big(\sum_{k=-\infty }^{\infty
}2^{k\alpha p}\text{ }\left\Vert f\chi _{k}\right\Vert _{q}^{p}\Big)^{1/p},
\end{equation*}%
\textit{with the usual modifications made when }$p=\infty $\textit{\ and/or }%
$q=\infty $\textit{.}
\end{definition}

The spaces $\dot{K}_{q}^{\alpha ,p}$ are quasi-Banach spaces and if $\min
(p,q)\geq 1$ then $\dot{K}_{q}^{\alpha ,p}$ are Banach spaces. When $\alpha
=0$ and $0<p=q\leq \infty $ then $\dot{K}_{p}^{0,p}$ coincides with the
Lebesgue spaces $L^{p}$. A detailed discussion of the properties of these
spaces my be found in the papers \cite{HerYang99}, \cite{LuYang1.95}, \cite%
{LuYang2.95}, and references therein.

\ \ \ Now, we\ present the Fourier analytical definition of Herz-type Besov
and Triebel-Lizorkin spaces\ and recall their basic properties. We first
need the concept of a smooth dyadic resolution of unity. Let $\phi _{0}$\ be
a function\ in $\mathcal{S}(\mathbb{R}^{n})$\ satisfying $\phi _{0}(x)=1$\
for\ $\left\vert x\right\vert \leq 1$\ and\ $\phi _{0}(x)=0$\ for\ $%
\left\vert x\right\vert \geq 2$.\ We put $\phi _{j}(x)=\phi
_{0}(2^{-j}x)-\phi _{0}(2^{1-j}x)$ for $j=1,2,3,...$. Then $\{\phi
_{j}\}_{j\in \mathbb{N}_{0}}$\ is a resolution of unity, $\sum_{j=0}^{\infty
}\phi _{j}(x)=1$ for all $x\in \mathbb{R}^{n}$.\ Thus we obtain the
Littlewood-Paley decomposition $f=\sum_{j=0}^{\infty }\mathcal{F}^{-1}\phi
_{j}\ast f$ of all $f\in \mathcal{S}^{\prime }(\mathbb{R}^{n})$ $($%
convergence in $\mathcal{S}^{\prime }(\mathbb{R}^{n}))$.

We are now in a position to state the definition of Herz-type Besov and
Triebel-Lizorkin spaces.

\begin{definition}
\label{Herz-Besov-Triebel}\textit{Let }$\alpha ,s\in \mathbb{R},0<p,q\leq
\infty $\textit{\ and }$0<\beta \leq \infty $\textit{. }

\begin{enumerate}
\item[$\mathrm{(i)}$] The \textit{Herz-type }Besov space $\dot{K}%
_{q}^{\alpha ,p}B_{\beta }^{s}$\ is the collection of all $f\in \mathcal{S}%
^{\prime }(\mathbb{R}^{n})$\ such that 
\begin{equation*}
\big\|f\big\|_{\dot{K}_{q}^{\alpha ,p}B_{\beta }^{s}}=\Big(%
\sum\limits_{j=0}^{\infty }2^{js\beta }\big\|\mathcal{F}^{-1}\phi _{j}\ast f%
\big\|_{\dot{K}_{q}^{\alpha ,p}}^{\beta }\Big)^{1/\beta }<\infty ,
\end{equation*}%
\textit{with the obvious modification if }$\beta =\infty .$

\item[$\mathrm{(ii)}$] Let $0<p,q<\infty $. \textit{The Herz-type
Triebel-Lizorkin space }$\dot{K}_{q}^{\alpha ,p}F_{\beta }^{s}$ \textit{is
the collection of all} $f\in \mathcal{S}^{\prime }(\mathbb{R}^{n})$\textit{\
such that}%
\begin{equation}
\big\|f\big\|_{\dot{K}_{q}^{\alpha ,p}F_{\beta }^{s}}=\Big\|\Big(%
\sum\limits_{j=0}^{\infty }2^{js\beta }\left\vert \mathcal{F}^{-1}\phi
_{j}\ast f\right\vert ^{\beta }\Big)^{1/\beta }\Big\|_{\dot{K}_{q}^{\alpha
,p}}<\infty ,  \label{H-T-L-space}
\end{equation}%
\textit{with the obvious modification if }$\beta =\infty .$
\end{enumerate}
\end{definition}

\begin{remark}
$\mathrm{Let}$\textrm{\ }$s\in \mathbb{R},0<p,q\leq \infty ,0<\beta \leq
\infty $ $\mathrm{and}$ $\alpha >-n/q$\textrm{. }$\mathrm{The}$ $\mathrm{%
spaces}$\textrm{\ }$\dot{K}_{q}^{\alpha ,p}B_{\beta }^{s}$ $\mathrm{and}$ $%
\dot{K}_{q}^{\alpha ,p}F_{\beta }^{s}$\textrm{\ }$\mathrm{are}$ $\mathrm{%
independent}$ $\mathrm{of}$ $\mathrm{the}$ $\mathrm{particular}$ $\mathrm{%
choice}$ $\mathrm{of}$ $\mathrm{the}$ $\mathrm{smooth}$ $\mathrm{dyadic}$ $%
\mathrm{resolution}$ $\mathrm{of}$ $\mathrm{unity}$\textrm{\ }$\{\phi
_{j}\}_{j\in \mathbb{N}_{0}}$ $\mathrm{(in}$ $\mathrm{the}$ $\mathrm{sense}$ 
$\mathrm{of\ equivalent}$ $\mathrm{quasi}$\textrm{-}$\mathrm{norms)}$\textrm{%
. }$\mathrm{In}$ $\mathrm{particular}$ $\dot{K}_{q}^{\alpha ,p}B_{\beta
}^{s} $ $\mathrm{and}$ $\dot{K}_{q}^{\alpha ,p}F_{\beta }^{s}$\textrm{\ }$%
\mathrm{are}$ $\mathrm{quasi}$\textrm{-}$\mathrm{Banach}$ $\mathrm{spaces}$ $%
\mathrm{and}$ $\mathrm{if}$\textrm{\ }$p,q,\beta \geq 1\mathrm{,}$ $\mathrm{%
then}$ $\dot{K}_{q}^{\alpha ,p}B_{\beta }^{s}$ $\mathrm{and}$ $\dot{K}%
_{q}^{\alpha ,p}F_{\beta }^{s}$\ $\mathrm{are}$ $\mathrm{Banach}$ $\mathrm{%
spaces}$\textrm{. }$\mathrm{Further}$ $\mathrm{results,concerning,}$ $%
\mathrm{for}$ $\mathrm{instance,}$ $\mathrm{lifting}$ $\mathrm{properties,}$ 
$\mathrm{Fourier}$ $\mathrm{multiplier}$ $\mathrm{and}$ $\mathrm{local}$ $%
\mathrm{means}$ $\mathrm{characterizations}$ $\mathrm{can}$ $\mathrm{be}$ $%
\mathrm{found}$ $\mathrm{in}$ $\mathrm{\cite{XuYang05},}$ $\mathrm{\cite%
{Xu03},}$ $\mathrm{\cite{Xu04}}$ $\mathrm{and}$ $\mathrm{\cite{Xu09}.}$
\end{remark}

Now we give the definitions of the spaces $B_{p,\beta }^{s}$ and $F_{p,\beta
}^{s}$.\vskip5pt

\begin{definition}
$\mathrm{(i)}$\textit{\ Let }$s\in \mathbb{R}$\textit{\ and }$0<p,\beta \leq
\infty $\textit{. The Besov space }$B_{p,\beta }^{s}$ \textit{is the
collection of all} $f\in \mathcal{S}^{\prime }(\mathbb{R}^{n})$\textit{\
such that} 
\begin{equation*}
\big\|f\big\|_{B_{p,\beta }^{s}}=\Big(\sum\limits_{j=0}^{\infty }2^{js\beta
}\left\Vert \mathcal{F}^{-1}\phi _{j}\ast f\right\Vert _{p}^{\beta }\Big)%
^{1/\beta }<\infty .
\end{equation*}%
$\mathrm{(ii)}$\textit{\ Let }$s\in \mathbb{R},0<p<\infty $\textit{\ and }$%
0<\beta \leq \infty $\textit{. The Triebel-Lizorkin space }$F_{p,\beta }^{s}$
\textit{is the collection of all} $f\in \mathcal{S}^{\prime }(\mathbb{R}%
^{n}) $\textit{\ such that}%
\begin{equation*}
\big\|f\big\|_{F_{p,\beta }^{s}}=\Big\|\Big(\sum\limits_{j=0}^{\infty
}2^{js\beta }\left\vert \mathcal{F}^{-1}\phi _{j}\ast f\right\vert ^{\beta }%
\Big)^{1/\beta }\Big\|_{p}<\infty .
\end{equation*}
\end{definition}

The theory of the spaces $B_{p,\beta }^{s}$ and $F_{p,\beta }^{s}$ has been
developed in detail in \cite{Triebel83}, \cite{Triebel92} and \cite%
{Triebel06} but has a longer history already including many contributors; we
do not want to discuss this here. Clearly, for $s\in \mathbb{R},0<p<\infty $
and $0<\beta \leq \infty ,$%
\begin{equation*}
\dot{K}_{p}^{0,p}F_{\beta }^{s}=F_{p,\beta }^{s}.
\end{equation*}

We introduce the sequence spaces associated with the function spaces $\dot{K}%
_{q}^{\alpha ,p}F_{\beta }^{s}$. If 
\begin{equation*}
\lambda =\{\lambda _{v,m}\in \mathbb{C}:v\in \mathbb{N}_{0},m\in \mathbb{Z}%
^{n}\},
\end{equation*}%
$\alpha ,s\in \mathbb{R},0<p,q\leq \infty $\ and $0<\beta \leq \infty $, we
set%
\begin{equation*}
\big\|\lambda \big\|_{\dot{K}_{q}^{\alpha ,p}b_{\beta }^{s}}=\Big(%
\sum_{v=0}^{\infty }2^{vs\beta }\Big\|\sum\limits_{m\in \mathbb{Z}%
^{n}}\lambda _{v,m}\chi _{v,m}\Big\|_{\dot{K}_{q}^{\alpha ,p}}^{\beta }\Big)%
^{1/\beta }
\end{equation*}%
and, with $0<p,q<\infty $, 
\begin{equation}
\big\|\lambda \big\|_{\dot{K}_{q}^{\alpha ,p}f_{\beta }^{s}}=\Big\|\Big(%
\sum_{v=0}^{\infty }\sum\limits_{m\in \mathbb{Z}^{n}}2^{vs\beta }|\lambda
_{v,m}|^{\beta }\chi _{v,m}\Big)^{1/\beta }\Big\|_{\dot{K}_{q}^{\alpha ,p}}.
\label{atomic-norm}
\end{equation}%
Let $\Phi ,\psi ,\varphi $ and $\Psi $ satisfy 
\begin{equation}
\Phi ,\Psi ,\varphi ,\psi \in \mathcal{S}(\mathbb{R}^{n})  \label{Ass1}
\end{equation}%
\begin{equation}
\text{supp}\mathcal{F}\Phi \text{, supp}\mathcal{F}\Psi \subset \overline{%
B(0,2)},\text{ \ }|\mathcal{F}\Phi (\xi )|,|\mathcal{F}\Psi (\xi )|\geq c%
\text{ if }|\xi |\leq \frac{5}{3},  \label{Ass2}
\end{equation}%
and 
\begin{equation}
\text{supp}\mathcal{F}\varphi \text{, supp}\mathcal{F}\psi \subset \overline{%
B(0,2)}\backslash B(0,1/2),\text{ \ }|\mathcal{F}\varphi (\xi )|,|\mathcal{F}%
\psi (\xi )|\geq c\text{ if }\frac{3}{5}\leq |\xi |\leq \frac{5}{3}
\label{Ass3}
\end{equation}%
such that%
\begin{equation}
\overline{\mathcal{F}\Phi (\xi )}\mathcal{F}\Psi (\xi )+\sum_{j=1}^{\infty }%
\overline{\mathcal{F}\varphi (2^{-j}\xi )}\mathcal{F}\psi (2^{-j}\xi
)=1,\quad \xi \in \mathbb{R}^{n},  \label{Ass4}
\end{equation}%
where $c>0$. Recall that the $\varphi $-transform $S_{\varphi }$ is defined
by setting $(S_{\varphi }f)_{0,m}=\langle f,\Psi _{m}\rangle $ where $\Psi
_{m}(x)=\Psi (x-m)$ and $(S_{\varphi }f)_{v,m}=\langle f,\varphi
_{v,m}\rangle $ where $\varphi _{v,m}(x)=2^{vn/2}\varphi (2^{v}x-m)$ and $%
v\in \mathbb{N}$. The inverse $\varphi $-transform $T_{\psi }$ is defined by 
\begin{equation*}
T_{\psi }\lambda =\sum_{m\in \mathbb{Z}^{n}}\lambda _{0,m}\Psi
_{m}+\sum_{v=1}^{\infty }\sum_{m\in \mathbb{Z}^{n}}\lambda _{v,m}\psi _{v,m},
\end{equation*}%
where $\lambda =\{\lambda _{v,m}\in \mathbb{C}:v\in \mathbb{N}_{0},m\in 
\mathbb{Z}^{n}\}$, see \cite{FJ90}.

For simplicity, in what follows, we use $\dot{K}_{p}^{\alpha ,q}A_{\beta
}^{s}$ to denote either $\dot{K}_{p}^{\alpha ,q}F_{\beta }^{s}$ or $\dot{K}%
_{p}^{\alpha ,q}B_{\beta }^{s}$. If $\dot{K}_{p}^{\alpha ,q}A_{\beta }^{s}$
means $\dot{K}_{p}^{\alpha ,q}F_{\beta }^{s}$ then the case $p=\infty $ is
excluded.\ To prove the main results of this paper we need the following
theorem, see \cite{Drihem16}.

\begin{theorem}
\label{phi-tran}\textit{Let }$\alpha ,s\in \mathbb{R},0<p,q\leq \infty
,0<\beta \leq \infty $ \textit{and }$\alpha >-n/q$. \textit{Suppose that }$%
\varphi $ and $\Phi $ satisfy \eqref{Ass1}-\eqref{Ass4}. The operators $%
S_{\varphi }:\dot{K}_{q}^{\alpha ,p}A_{\beta }^{s}\rightarrow \dot{K}%
_{q}^{\alpha ,p}a_{\beta }^{s}$ and $T_{\psi }:\dot{K}_{q}^{\alpha
,p}a_{\beta }^{s}\rightarrow \dot{K}_{q}^{\alpha ,p}A_{\beta }^{s}$ are
bounded. Furthermore, $T_{\psi }\circ S_{\varphi }$ is the identity on $\dot{%
K}_{q}^{\alpha ,p}A_{\beta }^{s}$.
\end{theorem}

We end this section with one more lemma, which is basically a consequence of
Hardy's inequality in the sequence Lebesgue space $\ell _{q}$.

\begin{lemma}
\label{lq-inequality}\textit{Let }$0<a<1$\textit{\ and }$0<q\leq \infty $%
\textit{. Let }$\left\{ \varepsilon _{k}\right\} $\textit{\ be a sequences
of positive real numbers and denote }$\delta
_{k}=\sum_{j=0}^{k}a^{k-j}\varepsilon _{j}$ and $\eta
_{k}=\sum_{j=k}^{\infty }a^{j-k}\varepsilon _{j},k\in \mathbb{N}_{0}$. Then
there exists constant $c>0\ $\textit{depending only on }$a$\textit{\ and }$q$
such that%
\begin{equation*}
\Big(\sum\limits_{k=0}^{\infty }\delta _{k}^{q}\Big)^{1/q}+\Big(%
\sum\limits_{k=0}^{\infty }\eta _{k}^{q}\Big)^{1/q}\leq c\text{ }\Big(%
\sum\limits_{k=0}^{\infty }\varepsilon _{k}^{q}\Big)^{1/q}.
\end{equation*}
\end{lemma}

\section{Jawerth embedding}

The classical Jawerth embedding says that:%
\begin{equation*}
F_{q,\infty }^{s_{2}}\hookrightarrow B_{s,q}^{s_{1}}
\end{equation*}%
if $s_{1}-n/s=s_{2}-n/q$ and $0<q<s<\infty $, see e.g. \cite{Ja77}. We will
extend this embeddings to Herz-type Besov Triebel-Lizorkin spaces.\ We
follow some ideas of Vyb\'{\i}ral, \cite{Vybiral08}, where use the technique
of non-increasing rearrangement.

\begin{definition}
Let $\mu $ be the Lebesgue measure in $\mathbb{R}^{n}$. If $f$ is a
measurable function on $\mathbb{R}^{n}$, we define the non-increasing
rearrangement of $f$ through%
\begin{equation*}
f^{\ast }(t)=\sup \{\lambda >0:m_{f}(\lambda )>t\}
\end{equation*}%
where $m_{f}$ is the distribution function of $f$.
\end{definition}

We shall use the following properties. If $0<p<\infty $, then%
\begin{equation}
\big\|f\big\|_{p}=\big\|f^{\ast }\mid L^{p}(0,\infty )\big\|
\label{property1}
\end{equation}%
for every measurable function $f$. Let $f$ and $g$ be two non-negative
measurable functions on $\mathbb{R}^{n}$. If $1\leq p\leq \infty $, then%
\begin{equation}
\big\|f+g\big\|_{p}\leq \big\|f^{\ast }+g^{\ast }\mid L^{p}(0,\infty )\big\|.
\label{property2}
\end{equation}%
The proof follows from Theorems 3.4 and 4.6 in \cite{BS88}. First, we will
prove the discrete version of Jawerth embedding.

\begin{theorem}
\label{embeddings6}\textit{Let }$\alpha _{1},\alpha _{2},s_{1},s_{2}\in 
\mathbb{R}$, $0<s,p\leq \infty ,0<q,r<\infty ,\alpha _{1}>-n/s\ $\textit{and 
}$\alpha _{2}>-n/q$. \textit{We suppose that }%
\begin{equation}
s_{1}-n/s-\alpha _{1}=s_{2}-n/q-\alpha _{2}.  \label{new-exp1}
\end{equation}

\noindent \textit{Under the following assumptions} \textit{\ }%
\begin{equation}
0<q<s\leq \infty ,q\leq r\text{\ and }\alpha _{2}>\alpha _{1}
\label{new-exp1.1}
\end{equation}%
or%
\begin{equation}
0<q<\min (s,p),q\leq r\leq \min (s,p)\ \text{and }\alpha _{2}=\alpha _{1}
\label{new-exp1.2}
\end{equation}%
or%
\begin{equation}
0<s\leq q<\infty ,\alpha _{2}+n/q>\alpha _{1}+n/s  \label{new-exp2}
\end{equation}%
or%
\begin{equation*}
0<s\leq q<\infty \text{, }q\leq r\leq p\leq \infty \text{\ and }\alpha
_{2}+n/q=\alpha _{1}+n/s,
\end{equation*}%
we have%
\begin{equation}
\dot{K}_{q}^{\alpha _{2},r}f_{\theta }^{s_{2}}\hookrightarrow \dot{K}%
_{s}^{\alpha _{1},p}b_{r}^{s_{1}},  \label{FJ-emb1}
\end{equation}%
where%
\begin{equation*}
\theta =\left\{ 
\begin{array}{ccc}
r & \text{if } & 0<s\leq q<\infty \text{, }q\leq r\leq p\leq \infty \text{
and }\alpha _{2}+n/q=\alpha _{1}+n/s \\ 
\infty &  & \text{otherwise.}%
\end{array}%
\right.
\end{equation*}
\end{theorem}

\textbf{Proof.} Let $\lambda \in \dot{K}_{q}^{\alpha _{2},r}f_{\theta
}^{s_{2}}$. We have%
\begin{eqnarray*}
\big\|\lambda \big\|_{\dot{K}_{s}^{\alpha _{1},p}b_{r}^{s_{1}}}^{r}
&=&\sum_{v=0}^{\infty }\Big(\sum\limits_{k=-\infty }^{\infty }2^{(k\alpha
_{1}+vs_{1})p}\big\|\sum\limits_{m\in \mathbb{Z}^{n}}\lambda _{v,m}\chi
_{v,m}\chi _{k}\big\|_{s}^{p}\Big)^{r/p} \\
&\leq &\sum_{v=0}^{\infty }\Big(\sum\limits_{k=-\infty }^{-v}\cdot \cdot
\cdot \Big)^{r/p}+\sum_{v=0}^{\infty }\Big(\sum\limits_{k=-v}^{\infty }\cdot
\cdot \cdot \Big)^{r/p} \\
&=&I+II.
\end{eqnarray*}

\textbf{Step 1. }We prove our embedding under the assumption %
\eqref{new-exp1.1}\ and we will estimate $I$ and $II$, respectively. We will
treat the case only where $0<s,p<\infty $. The case $s=\infty $ follows by
the embedding 
\begin{equation}
\dot{K}_{p_{0}}^{\alpha _{1},p}b_{r}^{s_{1}+\frac{n}{p_{0}}-\frac{n}{s}%
}\hookrightarrow \dot{K}_{s}^{\alpha _{1},p}b_{r}^{s_{1}}  \label{our-emb}
\end{equation}%
for some $0<q<p_{0}<s\leq \infty $, see \cite[Theorem\ 5.9]{Drihem13}.

\textbf{Estimation of }$I$\textit{.} Let $x\in C_{k}\cap Q_{v,m}$ and $y\in
Q_{v,m}$. We have $|x-y|\leq 2\sqrt{n}2^{-v}<2^{c_{n}-v}$ and from this it
follows that $|y|<2^{c_{n}-v}+2^{k}\leq 2^{c_{n}-v+2}$, which implies that $%
y $ is located in some ball $B(0,2^{c_{n}-v+2})$ and 
\begin{equation*}
|\lambda _{v,m}|^{t}\lesssim 2^{nv}\int_{B(0,2^{c_{n}-v+2})}|\lambda
_{v,m}|^{t}\chi _{v,m}(y)dy,
\end{equation*}%
where $t>0$. Then for any $x\in C_{k}$ we obtain%
\begin{eqnarray*}
\sum_{m\in \mathbb{Z}^{n}}|\lambda _{v,m}|^{t}\chi _{v,m}(x) &\lesssim
&2^{nv}\int_{B(0,2^{c_{n}-v+2})}\sum_{m\in \mathbb{Z}^{n}}|\lambda
_{v,m}|^{t}\chi _{v,m}(y)dy \\
&=&2^{nv}\big\|\sum_{m\in \mathbb{Z}^{n}}\lambda _{v,m}\chi _{v,m}\chi
_{B(0,2^{c_{n}-v+2})}\big\|_{t}^{t}.
\end{eqnarray*}%
Consequently,%
\begin{equation*}
2^{k\alpha _{1}+vs_{1}}\Big\|\sum_{m\in \mathbb{Z}^{n}}\lambda _{v,m}\chi
_{v,m}\chi _{k}\Big\|_{s}\lesssim \text{ }2^{v(s_{1}+\frac{n}{t})+k(\alpha
_{1}+\frac{n}{s})}\Big\|\sum_{m\in \mathbb{Z}^{n}}\lambda _{v,m}\chi
_{v,m}\chi _{B(0,2^{c_{n}-v+2})}\Big\|_{t}.
\end{equation*}%
We may choose $t>0$ such that $\frac{1}{t}>\max (\frac{1}{q},\frac{1}{q}+%
\frac{\alpha _{2}}{n})$. Therefore, since $\alpha _{1}+\frac{n}{s}>0$, 
\begin{equation*}
I\lesssim \sum_{v=0}^{\infty }2^{v(s_{1}+\frac{n}{t}-\alpha _{1}-\frac{n}{s}%
-s_{2})r}\sup_{j\geq 0}2^{s_{2}jr}\Big\|\sum_{m\in \mathbb{Z}^{n}}\lambda
_{j,m}\chi _{j,m}\chi _{B(0,2^{c_{n}-v+2})}\Big\|_{t}^{r},
\end{equation*}%
which can be estimated by, using \eqref{new-exp1},%
\begin{equation*}
c\sum_{v=0}^{\infty }2^{v\frac{nr}{d}}\Big(\sum_{i=-\infty }^{-v}2^{i\frac{%
n\sigma }{d}+\alpha _{2}\sigma i}\sup_{j\geq 0}2^{js_{2}\sigma }\Big\|%
\sum_{m\in \mathbb{Z}^{n}}\lambda _{j,m}\chi _{j,m}\chi _{i+c_{n}+2}\Big\|%
_{q}^{\sigma }\Big)^{r/\sigma },
\end{equation*}%
by H\"{o}lder's inequality, with $\sigma =\min (1,t)$ and $\frac{n}{d}=\frac{%
n}{t}-\frac{n}{q}-\alpha _{2}$. Hence Lemma \ref{lq-inequality} implies that 
\begin{equation*}
I\lesssim \sum_{i=0}^{\infty }2^{-\alpha _{2}ir}\sup_{j\geq 0}2^{js_{2}r}%
\Big\|\sum\limits_{m\in \mathbb{Z}^{n}}\lambda _{j,m}\chi _{j,m}\chi
_{2-i+c_{n}}\Big\|_{q}^{r}\lesssim \left\Vert \lambda \right\Vert _{\dot{K}%
_{q}^{\alpha _{2},r}f_{\infty }^{s_{2}}}^{r}.
\end{equation*}%
\textbf{Estimation of }$II$\textit{.} Our estimate use partially some
decomposition techniques already used in \cite{Vybiral08}. Set 
\begin{equation*}
h_{k}(x)=\sup_{v\geq 0}2^{vs_{2}}\sum\limits_{m\in \mathbb{Z}^{n}}\left\vert
\lambda _{v,m}\right\vert \chi _{v,m}(x)\chi _{k}(x).
\end{equation*}%
Then 
\begin{equation*}
\left\Vert \lambda \right\Vert _{\dot{K}_{q}^{\alpha _{2},r}f_{\infty
}^{s_{2}}}=\Big(\sum\limits_{k=-\infty }^{\infty }2^{k\alpha
_{2}r}\left\Vert h_{k}\right\Vert _{q}^{r}\Big)^{1/r}
\end{equation*}%
and 
\begin{equation*}
\left\vert \lambda _{v,m}\right\vert \leq 2^{-vs_{2}}\inf_{x\in
Q_{v,m}}h_{k}(x),\quad v\in \mathbb{N}_{0},m\in \mathbb{Z}^{n}.
\end{equation*}%
Using the fact that $\alpha _{2}>\alpha _{1}$ and the assumption %
\eqref{new-exp1} we estimate $II$ by%
\begin{eqnarray}
&&\sum_{v=0}^{\infty }2^{v(\frac{n}{s}-\frac{n}{q}+s_{2})r}\sup_{k\geq
-v}2^{k\alpha _{2}r}\Big\|\sum\limits_{m\in \mathbb{Z}^{n}}\lambda
_{v,m}\chi _{v,m}\chi _{k}\Big\|_{s}^{r}  \notag \\
&\leq &\sum\limits_{k=-\infty }^{\infty }2^{k\alpha _{2}r}\Big(%
\sum_{v=0}^{\infty }2^{v(\frac{n}{s}-\frac{n}{q}+s_{2})rt}\Big\|%
\sum\limits_{m\in \mathbb{Z}^{n}}\lambda _{v,m}\chi _{v,m}\chi _{k}\Big\|%
_{s}^{rt}\Big)^{1/t},  \label{sum}
\end{eqnarray}%
where 
\begin{equation*}
\frac{q}{r}<t<\min (1,\frac{s}{r})\text{ if }q<r\text{ and }t=1\text{ if }%
q=r.
\end{equation*}%
We can easly prove the estimate:%
\begin{equation*}
2^{vs_{2}s}\Big\|\sum\limits_{m\in \mathbb{Z}^{n}}\lambda _{v,m}\chi
_{v,m}\chi _{k}\Big\|_{s}^{s}\leq 2^{-vn}\sum\limits_{m\in \mathbb{Z}^{n}}%
\Big(\inf_{x\in Q_{v,m}}h_{k}(x)\Big)^{s}.
\end{equation*}%
Therefore, the sum $\sum_{v=0}^{\infty }\cdot \cdot \cdot $ in \eqref{sum}
can be estimated by 
\begin{eqnarray*}
&&\sum_{v=0}^{\infty }2^{-\frac{vnrt}{q}}\Big(\sum\limits_{m\in \mathbb{Z}%
^{n}}\Big(\inf_{x\in Q_{v,m}}h_{k}(x)\Big)^{s}\Big)^{rt/s} \\
&\leq &\sum_{v=0}^{\infty }2^{-\frac{vnrt}{q}}\Big(\sum\limits_{i=1}^{\infty
}\left( \left( h_{k}\right) ^{\ast }(2^{-vn}i)\right) ^{s}\Big)^{rt/s}.
\end{eqnarray*}%
Using the monotonicity of $h$ and the inequality $rt<s$, the last term is
bounded by%
\begin{eqnarray*}
&&c\sum_{v=0}^{\infty }2^{-\frac{vnrt}{q}}\Big(\sum\limits_{l=0}^{\infty
}2^{nl}\left( \left( h_{k}\right) ^{\ast }(2^{(l-v)n})\right) ^{s}\Big)%
^{rt/s} \\
&\lesssim &\sum_{v=0}^{\infty }2^{-\frac{vnrt}{q}}\sum\limits_{l=0}^{\infty
}2^{nl\frac{rt}{s}}\left( \left( h_{k}\right) ^{\ast }(2^{(l-v)n})\right)
^{rt} \\
&=&c\sum_{v=0}^{\infty }2^{-\frac{vnrt}{q}}\sum\limits_{j=-v}^{\infty
}2^{n(j+v)\frac{rt}{s}}\left( \left( h_{k}\right) ^{\ast }(2^{nj})\right)
^{rt} \\
&=&c\sum_{j=-\infty }^{\infty }2^{\frac{jnrt}{s}}\left( \left( h_{k}\right)
^{\ast }(2^{nj})\right) ^{r}\sum\limits_{v=-j}^{\infty }2^{nv(\frac{1}{s}-%
\frac{1}{q})rt} \\
&=&c\sum_{j=-\infty }^{\infty }2^{\frac{jnrt}{q}}\left( \left( h_{k}\right)
^{\ast }(2^{nj})\right) ^{rt}
\end{eqnarray*}%
Since $q<rt$, using the embedding $\ell _{q}\hookrightarrow \ell _{rt}$, we
get 
\begin{equation*}
\sum_{j=-\infty }^{\infty }2^{\frac{njrt}{q}}\left( \left( h_{k}\right)
^{\ast }(2^{nj})\right) ^{rt}\leq \Big(\sum_{j=-\infty }^{\infty
}2^{nj}\left( \left( h_{k}\right) ^{\ast }(2^{nj})\right) ^{q}\Big)%
^{rt/q}=\left\Vert h_{k}\right\Vert _{q}^{rt}.
\end{equation*}%
Consequently, we obtain $II\lesssim \left\Vert \lambda \right\Vert _{\dot{K}%
_{q}^{\alpha _{2},r}f_{\infty }^{s_{2}}}.$

\textbf{Step 2. }We prove our embedding under the assumption %
\eqref{new-exp1.2} and we need only to estimate $II$. Since $q\leq r\leq
\min (s,p)$, we can estimate $II$ by \eqref{sum} with $1$ in place of $t$.
Using similar arguments\ of Step 1, we get the desired estimate. Notice that
the case $s=\infty $ follows by the embedding \eqref{our-emb}\ for some $%
0<q<p_{0}<\min (s,p)\leq \infty $.

\noindent \textbf{Step 3. }We prove our embedding under the assumption %
\eqref{new-exp2} and again we need only to estimate $II$. By H\"{o}lder's
inequality we get%
\begin{equation*}
2^{vs_{1}}\Big\|\sum\limits_{m\in \mathbb{Z}^{n}}|\lambda _{v,m}|\chi
_{v,m}\chi _{k}\Big\|_{s}\leq 2^{(\frac{n}{s}-\frac{n}{q})k+vs_{1}}\Big\|%
\sum\limits_{m\in \mathbb{Z}^{n}}\lambda _{v,m}\chi _{v,m}\chi _{k}\Big\|%
_{q}.
\end{equation*}%
Hence $II$ can be\ estimated by%
\begin{eqnarray*}
&&c\sum_{v=0}^{\infty }2^{vs_{1}r}\Big(\sum\limits_{k=-v}^{\infty
}2^{k(\alpha _{1}+\frac{n}{s}-\frac{n}{q})p}\Big\|\sum\limits_{m\in \mathbb{Z%
}^{n}}\lambda _{v,m}\chi _{v,m}\chi _{k}\Big\|_{q}^{p}\Big)^{r/p} \\
&\leq &\sum_{v=0}^{\infty }2^{vs_{2}r}\Big(\sum\limits_{k=-v}^{\infty
}2^{(k+v)(\alpha _{1}-\alpha _{2}+\frac{n}{s}-\frac{n}{q})p}2^{k\alpha _{2}p}%
\Big\|\sum\limits_{m\in \mathbb{Z}^{n}}\lambda _{v,m}\chi _{v,m}\chi _{k}%
\Big\|_{q}^{p}\Big)^{r/p} \\
&\leq &\sum_{v=0}^{\infty }\Big(\sum\limits_{k=-v}^{\infty }2^{(k+v)(\alpha
_{1}-\alpha _{2}+\frac{n}{s}-\frac{n}{q})p}2^{k\alpha _{2}p}\Big\|%
\sup_{j\geq 0}2^{js_{2}}\sum\limits_{m\in \mathbb{Z}^{n}}\lambda _{j,m}\chi
_{j,m}\chi _{k}\Big\|_{q}^{p}\Big)^{r/p} \\
&\lesssim &\sum_{v=-\infty }^{\infty }2^{k\alpha _{2}r}\Big\|\sup_{j\geq
0}2^{js_{2}}\sum\limits_{m\in \mathbb{Z}^{n}}|\lambda _{j,m}|\chi _{j,m}\chi
_{k}\Big\|_{q}^{r} \\
&\lesssim &\left\Vert \lambda \right\Vert _{\dot{K}_{q}^{\alpha
_{2},r}f_{\infty }^{s_{2}}}^{r},
\end{eqnarray*}%
by Lemma\ \ref{lq-inequality}. If $\alpha _{2}+n/q=\alpha _{1}+n/s$\ and $%
r\leq p$, then 
\begin{eqnarray*}
II &\lesssim &\sum_{v=0}^{\infty }2^{vs_{2}r}\sum\limits_{k=-v}^{\infty
}2^{k\alpha _{2}r}\Big\|\sum\limits_{m\in \mathbb{Z}^{n}}\lambda _{v,m}\chi
_{v,m}\chi _{k}\Big\|_{q}^{r} \\
&\lesssim &\sum\limits_{k=-\infty }^{\infty }2^{k\alpha _{2}r}\Big\|\Big(%
\sum_{v=0}^{\infty }2^{vs_{2}r}\sum\limits_{m\in \mathbb{Z}^{n}}|\lambda
_{v,m}|^{r}\chi _{v,m}\chi _{k}\Big)^{1/r}\Big\|_{q}^{r} \\
&\leq &\left\Vert \lambda \right\Vert _{\dot{K}_{q}^{\alpha
_{2},r}f_{r}^{s_{2}}}^{r}.
\end{eqnarray*}%
The proof is complete. \ \ \rule{3mm}{3mm}

We would like to mention that\ $r$ on the right hand side of \eqref{FJ-emb1}
is optimal. Indeed, for $v\in \mathbb{N}_{0}$ and $N\geq 1$, we put 
\begin{equation*}
\lambda _{v,m}^{N}=\left\{ 
\begin{array}{ccc}
2^{-(s_{1}-\frac{1}{s}-\alpha _{1})v}\sum_{i=1}^{N}\chi _{i}(2^{v-1}) & 
\text{if} & m=1 \\ 
0 &  & \text{otherwise,}%
\end{array}%
\right.
\end{equation*}%
and $\lambda ^{N}=\{\lambda _{v,m}^{N}:v\in \mathbb{N}_{0},m\in \mathbb{Z}\}$%
. We have%
\begin{equation*}
\big\|\lambda ^{N}\big\|_{\dot{K}_{q}^{\alpha _{2},r}f_{\theta
}^{s_{2}}}^{r}=\sum_{k=-\infty }^{\infty }2^{\alpha _{2}kr}\Big\|\Big(%
\sum_{v=0}^{\infty }2^{vs_{2}\theta }|\lambda _{v,1}^{N}|^{\theta }\chi
_{v,1}\Big)^{1/\theta }\chi _{k}\Big\|_{q}^{r}.
\end{equation*}%
We can rewrite the last statement as follows: 
\begin{eqnarray*}
&&\sum_{k=1-N}^{0}2^{\alpha _{2}kr}\Big\|\Big(\sum_{v=1}^{N}2^{(s_{2}-s_{1}+%
\frac{1}{s}+\alpha _{1})v\theta }\chi _{v,1}\Big)^{1/\theta }\chi _{k}\Big\|%
_{q}^{r} \\
&=&\sum_{k=1-N}^{0}2^{\alpha _{2}kr}\Big\|2^{(s_{2}-s_{1}+\frac{1}{s}+\alpha
_{1})(1-k)}\chi _{1-k,1}\Big\|_{q}^{r} \\
&=&c\text{ }N,
\end{eqnarray*}%
where the constant $c>0$ does not depend on $N$. Now%
\begin{equation*}
\big\|\lambda ^{N}\big\|_{\dot{K}_{s}^{\alpha _{1},p}b_{\sigma
}^{s_{1}}}^{\sigma }=\sum_{v=0}^{\infty }2^{vs_{1}\sigma }\Big(%
\sum_{k=-\infty }^{\infty }2^{\alpha _{1}kp}\Big\|\sum\limits_{m\in \mathbb{Z%
}}|\lambda _{v,m}^{N}|\chi _{v,m}\chi _{k}\Big\|_{s}^{p}\Big)^{\sigma /p}.
\end{equation*}%
Again we can rewrite the last statement as follows:%
\begin{equation*}
\big\|\lambda ^{N}\big\|_{\dot{K}_{s}^{\alpha _{1},p}b_{\sigma
}^{s_{1}}}^{\sigma }=\sum_{v=1}^{N}2^{(\frac{1}{s}+\alpha _{1})v\sigma }\Big(%
\sum_{k=1-N}^{0}2^{\alpha _{1}kp}\Big\|\chi _{v,1}\chi _{k}\Big\|_{s}^{p}%
\Big)^{\sigma /p}=cN,
\end{equation*}%
where the constant $c>0$ does not depend on $N$. If the embeddings %
\eqref{FJ-emb1} holds then for any $N\in \mathbb{N}$, $N^{\frac{1}{\sigma }-%
\frac{1}{r}}\leq C$. Thus, we conclude that $0<r\leq \sigma <\infty $ must
necessarily hold by letting $N\rightarrow +\infty $.

Using Theorems \ref{phi-tran} and \ref{embeddings6}, we have the following
Jawerth embedding.

\begin{theorem}
\label{embeddings6.1}\textit{Let }$\alpha _{1},\alpha _{2},s_{1},s_{2}\in 
\mathbb{R}$, $0<s,p<\infty ,0<q,r\leq \infty ,\alpha _{1}>-n/s\ $\textit{and 
}$\alpha _{2}>-n/q$. Under the hypothesis of Therorem \ref{embeddings6}\ we
have%
\begin{equation}
\dot{K}_{q}^{\alpha _{2},r}F_{\theta }^{s_{2}}\hookrightarrow \dot{K}%
_{s}^{\alpha _{1},p}B_{r}^{s_{1}},  \label{Jawerth}
\end{equation}%
where%
\begin{equation*}
\theta =\left\{ 
\begin{array}{ccc}
r & \text{if } & 0<s\leq q<\infty \text{, }q\leq r\leq p\leq \infty \text{
and }\alpha _{2}+n/q=\alpha _{1}+n/s \\ 
\infty &  & \text{otherwise.}%
\end{array}%
\right.
\end{equation*}
\end{theorem}

From this theorem and the fact that $\dot{K}_{q}^{0,q}F_{\infty
}^{s_{2}}=F_{q,\infty }^{s_{2}}$ and $\dot{K}_{s}^{0,p}B_{q}^{s_{1}}%
\hookrightarrow \dot{K}_{s}^{0,s}B_{q}^{s_{1}}=B_{s,q}^{s_{2}}$, with $%
0<p<s<\infty $, we obtain the following embeddings 
\begin{equation*}
F_{q,\infty }^{s_{2}}\hookrightarrow \dot{K}_{s}^{0,p}B_{q}^{s_{1}}%
\hookrightarrow B_{s,q}^{s_{2}},
\end{equation*}%
if 
\begin{equation*}
0<q<p<s<\infty \text{ and }s_{1}-n/s=s_{2}-n/q.
\end{equation*}%
From this theorem and the fact that $\dot{K}_{q}^{\alpha ,r}F_{2}^{0}=\dot{K}%
_{q}^{\alpha ,r}$ for $1<r,q<\infty $ and $-\frac{n}{q}<\alpha <n-\frac{n}{q}
$, see \cite{XuYang03} and again $\dot{K}%
_{s}^{0,s}B_{r}^{s_{1}}=B_{s,r}^{s_{1}}$\ we immediately arrive at the
following embedding between Herz and Besov spaces.

\begin{theorem}
\textit{Let }$\alpha ,s_{1}\in \mathbb{R}$, $0<s\leq \infty ,1<r,q<\infty \ $%
\textit{and }$0\leq \alpha <n-\frac{n}{q}$. \textit{We suppose that }$\tfrac{%
n}{s}-s_{1}=\alpha +\tfrac{n}{q}$.\textit{\ Let}%
\begin{equation*}
0<q<s\leq \infty ,q\leq r\ \text{and }\alpha >0
\end{equation*}%
or%
\begin{equation*}
0<q<s\leq \infty ,q\leq r\leq s\ \text{and }\alpha =0
\end{equation*}%
or%
\begin{equation*}
0<s\leq q<\infty ,\alpha >\tfrac{n}{s}-\tfrac{n}{q}
\end{equation*}%
or%
\begin{equation*}
0<s\leq q\leq \infty \text{, }q\leq r\leq s\leq \infty \text{\ and }\alpha =%
\tfrac{n}{s}-\tfrac{n}{q}.
\end{equation*}%
Then%
\begin{equation*}
\dot{K}_{q}^{\alpha ,r}\hookrightarrow B_{s,r}^{s_{1}},
\end{equation*}%
where%
\begin{equation*}
r=2\text{ if }0<s\leq q<\infty \text{, }q\leq 2\leq s\leq \infty \text{ and }%
\alpha =\tfrac{n}{s}-\tfrac{n}{q}.
\end{equation*}
\end{theorem}

Some embeddings between Herz spaces and homogenous Besov spaces can be found
in \cite{T11}.

\section{Franke embedding}

The classical Franke embedding may be rewritten us follows: 
\begin{equation*}
B_{q,s}^{s_{2}}\hookrightarrow F_{s,\infty }^{s_{1}},
\end{equation*}%
if $s_{1}-n/s=s_{2}-n/q$ and $0<q<s<\infty $, see e.g. \cite{Fr86}. As in
Section 3 we will extend this embeddings to Herz-type Besov-Triebel-Lizorkin
spaces. Again, we follow some ideas of Vyb\'{\i}ral, \cite{Vybiral08}. We
will prove the discrete version of Franke embedding.

\begin{theorem}
\label{F-emb3}\textit{Let }$\alpha _{1},\alpha _{2},s_{1},s_{2}\in \mathbb{R}
$, $0<s,p,q<\infty ,0<\theta \leq \infty ,\alpha _{1}>-\frac{n}{s}\ $\textit{%
and }$\alpha _{2}>-\frac{n}{q}$. \textit{We suppose that }%
\begin{equation*}
s_{1}-\tfrac{n}{s}-\alpha _{1}=s_{2}-\tfrac{n}{q}-\alpha _{2}.
\end{equation*}

\noindent \textit{Let }%
\begin{equation}
0<q<s<\infty ,\alpha _{2}\geq \alpha _{1},  \label{Newexp21}
\end{equation}%
or 
\begin{equation}
0<s\leq q<\infty \text{\ and }\alpha _{2}+\tfrac{n}{q}>\alpha _{1}+\tfrac{n}{%
s}.  \label{newexp2}
\end{equation}%
Then%
\begin{equation}
\dot{K}_{q}^{\alpha _{2},p}b_{p}^{s_{2}}\hookrightarrow \dot{K}_{s}^{\alpha
_{1},p}f_{\theta }^{s_{1}}.  \label{Sobolev-emb1}
\end{equation}
\end{theorem}

\textbf{Proof}. We prove our embedding under the conditions \eqref{Newexp21}%
. Let $\lambda \in \dot{K}_{q}^{\alpha _{2},p}b_{p}^{s_{2}}$. We have%
\begin{eqnarray*}
\big\|\lambda \big\|_{\dot{K}_{s}^{\alpha _{1},p}f_{\theta }^{s_{1}}}^{p}
&=&\sum\limits_{k=-\infty }^{\infty }2^{k\alpha _{1}p}\Big\|\Big(%
\sum_{v=0}^{\infty }\sum\limits_{m\in \mathbb{Z}^{n}}2^{vs_{1}\theta
}|\lambda _{v,m}|^{\theta }\chi _{v,m}\chi _{k}\Big)^{1/\theta }\Big\|%
_{s}^{p} \\
&=&\sum\limits_{k=-\infty }^{0}\cdot \cdot \cdot +\sum\limits_{k=1}^{\infty
}\cdot \cdot \cdot \\
&=&J_{1}+J_{2}.
\end{eqnarray*}%
\textbf{Estimation of }$J_{1}$. Let $c_{n}=1+[\log _{2}(2\sqrt{n}+1)]$.
Obviously 
\begin{eqnarray*}
J_{1} &\lesssim &\sum\limits_{k=-\infty }^{0}2^{k\alpha _{1}p}\Big\|\Big(%
\sum_{v=0}^{c_{n}-k+1}\cdot \cdot \cdot \Big)^{1/\theta }\Big\|%
_{s}^{p}+\sum\limits_{k=-\infty }^{0}2^{k\alpha _{1}p}\Big\|\Big(%
\sum_{v=c_{n}-k+2}^{\infty }\cdot \cdot \cdot \Big)^{1/\theta }\Big\|_{s}^{p}
\\
&=&T_{1}+T_{2}.
\end{eqnarray*}%
The same analysis as in the proof of Theorem \ref{embeddings6} shows that 
\begin{equation*}
\sum_{m\in \mathbb{Z}^{n}}|\lambda _{v,m}|^{t}\chi _{v,m}(x)\lesssim 2^{nv}%
\Big\|\sum_{m\in \mathbb{Z}^{n}}|\lambda _{v,m}|\chi _{v,m}\chi
_{B(0,2^{c_{n}-v+2})}\Big\|_{t}^{t}
\end{equation*}%
for any $x\in C_{k}$. From Lemma \ref{lq-inequality}, since $\alpha _{1}+%
\frac{n}{s}>0$, $T_{1}$\ does not exceed 
\begin{equation*}
c\sum\limits_{v=0}^{\infty }2^{v(s_{1}-\alpha _{1}-\frac{n}{s}+\frac{n}{t})p}%
\Big\|\sum_{m\in \mathbb{Z}^{n}}|\lambda _{v,m}|\chi _{v,m}\chi
_{B(0,2^{c_{n}-v+2})}\Big\|_{t}^{p}.
\end{equation*}%
We may choose $t>0$ such that $\frac{1}{t}>\max (\frac{1}{q},\frac{1}{q}+%
\frac{\alpha _{2}}{n})$, $\sigma =\min (1,t)$ and $\frac{n}{d}=\frac{n}{t}-%
\frac{n}{q}-\alpha _{2}$. By H\"{o}lder's inequality, this term is bounded by%
\begin{equation*}
c\sum\limits_{v=0}^{\infty }2^{v\frac{n}{d}p}\Big(\sum_{i=-v}^{\infty }2^{i%
\frac{n\sigma }{d}+\alpha _{2}\sigma i}\sup_{j\geq 0}\Big\|\sum_{m\in 
\mathbb{Z}^{n}}2^{s_{2}j}|\lambda _{j,m}|\chi _{j,m}\chi _{i+c_{n}+2}\Big\|%
_{q}^{\sigma }\Big)^{p/\sigma }.
\end{equation*}%
Using again Lemma \ref{lq-inequality}, the last term is bounded by%
\begin{equation*}
c\sum_{i=0}^{\infty }2^{-\alpha _{2}ip}\sup_{j\geq 0}\Big\|\sum\limits_{m\in 
\mathbb{Z}^{n}}2^{s_{2}j}|\lambda _{j,m}|\chi _{j,m}\chi _{2-i+c_{n}}\Big\|%
_{q}^{p}\lesssim \left\Vert \lambda \right\Vert _{\dot{K}_{q}^{\alpha
_{2},p}b_{p}^{s_{2}}}^{p}.
\end{equation*}%
Now we estimate $T_{2}$. First let us consider $\alpha _{2}>\alpha _{1}$. We
have%
\begin{eqnarray*}
T_{2} &\lesssim &\sum\limits_{k=-\infty }^{0}2^{k\alpha _{2}p}\sup_{v\geq
c_{n}+2-k}2^{v(s_{2}-\frac{n}{q}+\frac{n}{s})p}\Big\|\sum\limits_{m\in 
\mathbb{Z}^{n}}\lambda _{v,m}\chi _{v,m}\chi _{k}\Big\|_{s}^{p} \\
&=&\sum\limits_{k=-\infty }^{0}2^{k\alpha _{2}p}\sup_{v\geq
c_{n}+2-k}2^{v(s_{2}-\frac{n}{q}+\frac{n}{s})p}\left\Vert h_{v,k}\right\Vert
_{s}^{p}.
\end{eqnarray*}%
Let us prove that%
\begin{equation*}
2^{v(\frac{n}{s}-\frac{n}{q})}\left\Vert h_{v,k}\right\Vert _{s}\lesssim %
\Big\|\sum\limits_{m\in \mathbb{Z}^{n}}\lambda _{v,m}\chi _{v,m}\chi _{%
\tilde{C}_{k}}\Big\|_{q}=\delta ,
\end{equation*}%
where $\tilde{C}_{k}=\cup _{i=-1}^{2}C_{k+i}$, wich equivalent to%
\begin{equation*}
\int 2^{v(\frac{n}{s}-\frac{n}{q})s}(h_{v,k}(x)\delta ^{-1})^{s}dx\lesssim 1.
\end{equation*}%
Let $x\in C_{k}\cap Q_{v,m}$ and $y\in Q_{v,m}$\ with $v\geq c_{n}-k+2$. We
have $|x-y|\leq 2\sqrt{n}2^{-v}<2^{c_{n}-v}$ and from this it follows that $%
2^{k-2}<|y|<2^{c_{n}-v}+2^{k}<2^{k+2}$, which implies that $y$ is located in 
$\widetilde{C}_{k}$. Therefore,%
\begin{equation*}
|\lambda _{v,m}|^{q}\chi _{v,m}(x)\lesssim 2^{nv}\int_{\tilde{C}%
_{k}}|\lambda _{v,m}|^{q}\chi _{v,m}(y)dy,
\end{equation*}%
if $x\in C_{k}\cap Q_{v,m}$. We have%
\begin{eqnarray*}
\int 2^{v(\frac{n}{s}-\frac{n}{q})s}(h_{v,k}(x)\delta ^{-1})^{s}dx &=&\int
\left( 2^{-v\frac{n}{q}}h_{v,k}(x)\delta ^{-1}\right) ^{s-q}\left(
h_{v,k}(x)\delta ^{-1}\right) ^{q}dx \\
&\lesssim &\int \left( h_{v,k}(x)\delta ^{-1}\right) ^{q}dx \\
&\lesssim &1.
\end{eqnarray*}%
Consequently,%
\begin{equation*}
T_{2}\lesssim \sum\limits_{k=-\infty }^{0}2^{k\alpha _{2}p}\sup_{v\geq
0}2^{vs_{2}p}\left\Vert h_{v,k}\right\Vert _{q}^{p}\leq \left\Vert \lambda
\right\Vert _{\dot{K}_{q}^{\alpha _{2},p}b_{p}^{s_{2}}}^{p}.
\end{equation*}%
Now let us consider $\alpha _{2}=\alpha _{1}$. We can suppose that $\theta
\leq q$\ and $p\leq s$ , since the opposite cases can be obtained by the
fact that\ $\ell _{q}\hookrightarrow \ell _{\theta }$ and/or $\ell
_{s}\hookrightarrow \ell _{p}$, repectively. Observe that%
\begin{equation*}
\left\vert 2^{-v}m\right\vert \leq \left\vert x-2^{-v}m\right\vert
+\left\vert x\right\vert \leq 2^{k}+\sqrt{n}2^{-v}\leq 2^{k+1}
\end{equation*}%
and 
\begin{equation*}
\left\vert 2^{-v}m\right\vert \geq \left\vert \left\vert
x-2^{-v}m\right\vert -\left\vert x\right\vert \right\vert \geq 2^{k-1}-\sqrt{%
n}2^{-v}\geq 2^{k-2}
\end{equation*}%
if $x\in C_{k}\cap Q_{v,m}$ and $v\geq c_{n}+2-k$. Hence $m$ is located in%
\begin{equation*}
\bar{C}_{k+v}=\{m\in \mathbb{Z}^{n}:2^{k+v-2}\leq \left\vert m\right\vert
\leq 2^{k+v+1}\}.
\end{equation*}%
Therefore\ $T_{2}$ can be estimated by%
\begin{equation*}
\sum\limits_{k=-\infty }^{0}2^{k\alpha _{1}p}\Big\|\Big(\sum_{v=c_{n}+2-k}^{%
\infty }\sum_{m\in \bar{C}_{k+v}}2^{vs_{1}\theta }|\lambda _{v,m}|^{\theta
}\chi _{v,m}\Big)^{1/\theta }\Big\|_{s}^{p}.
\end{equation*}%
Let%
\begin{equation*}
\tilde{\lambda}_{v,m_{1}}^{1,k}=\max_{m\in \bar{C}_{k+v}}\left\vert \lambda
_{v,m}\right\vert ,\quad m_{1}\in \mathbb{Z}^{n}
\end{equation*}%
and (decreassing rearrangement of $\{\lambda _{v,m}\}_{m\in \bar{C}_{k+v}}$) 
\begin{equation*}
\tilde{\lambda}_{v,m_{j}}^{j,k}=\max_{m^{i}\in \bar{C}_{k+v},i=1,...,j}%
\sum_{i=1}^{j}\left\vert \lambda _{v,m^{i}}\right\vert -\sum_{i=1}^{j-1}%
\tilde{\lambda}_{v,m_{i}}^{i,k},\quad m_{j}\in \mathbb{Z}^{n},j\geq 2.
\end{equation*}%
Then%
\begin{equation*}
\sum_{m\in \bar{C}_{k+v}}|\lambda _{v,m}|\chi _{v,m}=\sum_{i=1}^{L_{k+v}}%
\tilde{\lambda}_{v,m_{i}}^{i,k}\chi _{v,m_{i}}=f
\end{equation*}%
for some $L_{k+v}\in \mathbb{N}$. It is not difficult to see that%
\begin{equation*}
f^{\ast }(t)=\sum_{i=1}^{L_{k+v}}|\tilde{\lambda}_{v,m_{i}}^{i,k}|\tilde{\chi%
}_{[B_{i-1},B_{i})}(t),
\end{equation*}%
with 
\begin{equation*}
B_{i}=\sum_{j=1}^{i}\left\vert Q_{v,m_{j}}\right\vert =2^{-vn}i,\quad
i=1,...,L_{k+v}.
\end{equation*}%
Using the properties \eqref{property1} and \eqref{property2} we can estimate 
$\big\|\cdot \cdot \cdot \big\|_{s}^{p}$ by%
\begin{equation*}
c\Big\|\sum_{v=c_{n}+2-k}^{\infty }\sum_{i=1}^{L_{k+v}}2^{vs_{1}\theta }|%
\tilde{\lambda}_{v,m_{i}}^{i,k}|^{\theta }\tilde{\chi}_{v,i}\mid L^{s/\theta
}(0,\infty )\Big\|^{p/\theta },
\end{equation*}%
where $\tilde{\chi}_{v,i}$ is a characteristic function of the interval $%
(2^{-v}(i-1),2^{-v}i)$. By duality, the last norm may be rewritten as%
\begin{equation}
\sup \int_{0}^{\infty }\sum_{v=c_{n}+2-k}^{\infty
}\sum_{i=1}^{L_{k+v}}2^{vs_{1}\theta }|\tilde{\lambda}_{v,m_{i}}^{i,k}|^{%
\theta }\tilde{\chi}_{v,i}(x)g(x)dx  \label{Dual}
\end{equation}%
where the supremum is taken over all non-increasing non-negative measurable
functions $g$ with $\left\Vert g\mid L^{\beta }(0,\infty )\right\Vert \leq 1$
and $\beta $ is the conjugated index to $s/\theta $. Similarly, $\varrho $\
stands for the conjugated index to $q/\theta $. Let%
\begin{equation*}
g_{v,i}=\int_{0}^{\infty }\tilde{\chi}_{v,i}(x)g(x)dx.
\end{equation*}%
H\"{o}lder's inequality implies that%
\begin{eqnarray*}
&&\sum_{v=c_{n}+2-k}^{\infty }\sum_{i=1}^{L_{k+v}}2^{vs_{1}\theta }|\tilde{%
\lambda}_{v,m_{i}}^{i,k}|^{\theta }\tilde{\chi}_{v,i}g_{v,i} \\
&\leq &\sum_{v=c_{n}+2-k}^{\infty }\Big(\sum_{i=1}^{L_{k+v}}2^{vs_{2}q}|%
\tilde{\lambda}_{v,m_{i}}^{i,k}|^{q}\Big)^{\theta /q}\Big(%
\sum\limits_{h=1}^{\infty }2^{v(s_{1}-s_{2})\theta \varrho }g_{v,h}^{\varrho
}\Big)^{1/\varrho } \\
&\leq &\Big(\sum_{v=c_{n}+2-k}^{\infty }\Big(\sum_{i=1}^{L_{k+v}}2^{v(s_{2}-%
\frac{n}{q})q}|\tilde{\lambda}_{v,m_{i}}^{i,k}|^{q}\Big)^{s/q}\Big)^{\theta
/s}\Big(\sum_{v=c_{n}+2-k}^{\infty }\Big(\sum\limits_{h=1}^{\infty }2^{%
\tfrac{vn\theta \varrho }{s}}g_{v,h}^{\varrho }\Big)^{\beta /\varrho }\Big)%
^{1/\beta }.
\end{eqnarray*}%
As in \cite{Vybiral08} we can prove that the second term is bounded. Clearly
the first term can be estimated by%
\begin{equation*}
c\Big(\sum_{v=c_{n}+2-k}^{\infty }\Big(\sum_{m\in \bar{C}_{k+v}}2^{v(s_{2}-%
\frac{n}{q})q}|\lambda _{v,m}|^{q}\Big)^{s/q}\Big)^{\theta /s}\lesssim \Big(%
\sum_{v=1}^{\infty }2^{vs_{2}s}\Big\|\sum_{m\in \bar{C}_{k+v}}\lambda
_{v,m}\chi _{v,m}\chi _{\breve{C}_{k}}\Big\|_{q}^{s}\Big)^{\theta /s},
\end{equation*}%
where $\breve{C}_{k}=\cup _{i=-2}^{3}C_{k+i}$. Using the well-known
inequality%
\begin{equation}
\Big(\sum_{j=0}^{\infty }\left\vert a_{j}\right\vert \Big)^{\rho }\leq
\sum_{j=0}^{\infty }\left\vert a_{j}\right\vert ^{\rho },\quad \left\{
a_{j}\right\} _{j}\subset \mathbb{C},\text{ }\rho \in (0,1],  \label{lp-ine}
\end{equation}%
we obtain that $T_{2}$ can be estimated by $c\left\Vert \lambda \right\Vert
_{\dot{K}_{q}^{\alpha _{2},p}b_{p}^{s_{2}}}^{p}$.

\textbf{Estimation of }$J_{2}$. We use the same notations as in the\
estimation of $J_{1}$. We have%
\begin{equation*}
J_{2}\leq \sum\limits_{k=1}^{c_{n}+1}\cdot \cdot \cdot
+\sum\limits_{k=c_{n}+2}^{\infty }\cdot \cdot \cdot .
\end{equation*}%
As in the estimation of $T_{2}$, the second term can be estimated by $%
c\left\Vert \lambda \right\Vert _{\dot{K}_{q}^{\alpha
_{2},p}b_{p}^{s_{2}}}^{p}$. Now the first term is bounded by%
\begin{eqnarray*}
&&c\sum\limits_{k=1}^{c_{n}+1}2^{k\alpha _{1}p}\Big\|\Big(%
\sum_{v=0}^{c_{n}-k+1}\cdot \cdot \cdot \Big)^{1/\theta }\Big\|%
_{s}^{p}+\sum\limits_{k=1}^{c_{n}+1}2^{k\alpha _{1}p}\Big\|\Big(%
\sum_{v=c_{n}-k+2}^{\infty }\cdot \cdot \cdot \Big)^{1/\theta }\Big\|_{s}^{p}
\\
&\lesssim &\left\Vert \lambda \right\Vert _{\dot{K}_{q}^{\alpha
_{2},p}b_{p}^{s_{2}}}^{p},
\end{eqnarray*}%
where again we used the same arguments as in the estimation of $T_{1}$ and $%
T_{2}$.

Using a combination of the arguments used in Step 3 of the proof Theorem \ref%
{embeddings6} we prove our embedding under the conditions \eqref{newexp2}.
The proof is complete. \ \ \rule{3mm}{3mm}

Also as above $p$ on the right hand side of \eqref{Sobolev-emb1} is optimal.

Using Theorems \ref{phi-tran} and \ref{F-emb3}, we have the following Franke
embedding.

\begin{theorem}
\label{embeddings6.2}\textit{Let }$\alpha _{1},\alpha _{2},s_{1},s_{2}\in 
\mathbb{R}$, $0<s,q<\infty $, $0<p<\infty ,0<\theta \leq \infty ,\alpha
_{1}>-\frac{n}{s}\ $\textit{and }$\alpha _{2}>-\frac{n}{q}$. Under the
hypothesis of Therorem \ref{F-emb3}\ we have%
\begin{equation}
\dot{K}_{q}^{\alpha _{2},p}B_{p}^{s_{2}}\hookrightarrow \dot{K}_{s}^{\alpha
_{1},p}F_{\theta }^{s_{1}}.  \label{Frank}
\end{equation}
\end{theorem}

We would like to mention that from this theorem we have%
\begin{equation*}
B_{q,s}^{s_{2}}\hookrightarrow \dot{K}_{q}^{0,s}B_{s}^{s_{2}}\hookrightarrow
F_{s,\theta }^{s_{1}},
\end{equation*}%
if $0<q<s<\infty $, $0<\theta \leq \infty $ and\textit{\ }%
\begin{equation*}
s_{1}-\tfrac{n}{s}=s_{2}-\tfrac{n}{q}.
\end{equation*}
Also we immediately arrive at the following embedding between Herz and Besov
spaces.

\begin{theorem}
\textit{Let }$\alpha ,s_{2}\in \mathbb{R},1<s,q<\infty \ $\textit{and }$-%
\frac{n}{s}<\alpha \leq 0$. \textit{We suppose that }$\tfrac{n}{s}+\alpha =%
\tfrac{n}{q}-s_{2}$.\textit{\ Let }%
\begin{equation*}
1<q<s<\infty ,
\end{equation*}%
or 
\begin{equation*}
1<s\leq q<\infty \text{\ and }\alpha <\tfrac{n}{q}-\tfrac{n}{s}.
\end{equation*}%
Then%
\begin{equation*}
B_{q,q}^{s_{2}}\hookrightarrow \dot{K}_{s}^{\alpha ,q}.
\end{equation*}
\end{theorem}

By the same examples of \cite{Drihem13}, the assumptions $s_{1}-\tfrac{n}{s}%
-\alpha _{1}\leq s_{2}-\tfrac{n}{q}-\alpha _{2}\ $and $\alpha _{2}+\tfrac{n}{%
q}\geq \alpha _{1}+\tfrac{n}{s}$\ are necessary. Indeed, let $\eta \in 
\mathcal{S}\left( \mathbb{R}^{n}\right) $ be a function such that \textrm{%
supp}$\mathcal{F}\eta \subset \left\{ \xi \in \mathbb{R}^{n}:3/4<|\xi
|<1\right\} $. For $x\in \mathbb{R}^{n}$ and $N\in \mathbb{N}$ we put $%
f_{N}\left( x\right) =\eta \left( 2^{N}x\right) $. First we have $\eta \in 
\dot{K}_{s}^{\alpha _{1},p}\cap \dot{K}_{q}^{\alpha _{2},r}\cap \dot{K}%
_{q}^{\alpha _{2},p}$. Due to the support properties of the function $\eta $
we have for any $j\in \mathbb{N}_{0}$%
\begin{equation*}
\mathcal{F}^{-1}\phi _{j}\ast f_{N}=\left\{ 
\begin{array}{cc}
f_{N}, & j=N \\ 
0, & \text{otherwise}.%
\end{array}%
\right.
\end{equation*}%
Hence%
\begin{eqnarray*}
\left\Vert f_{N}\right\Vert _{\dot{K}_{s}^{\alpha _{1},p}B_{\beta }^{s_{1}}}
&=&2^{s_{1}N}\left\Vert f_{N}\right\Vert _{\dot{K}_{s}^{\alpha _{1},p}} \\
&=&2^{s_{1}N}\Big(\sum_{k=-\infty }^{\infty }2^{k\alpha _{1}p}\left\Vert
f_{N}\chi _{k}\right\Vert _{s}^{p}\Big)^{1/p} \\
&=&2^{\left( s_{1}-n/s\right) N}\Big(\sum_{k=-\infty }^{\infty }2^{k\alpha
_{1}p}\left\Vert \eta \chi _{k+N}\right\Vert _{s}^{p}\Big)^{1/p} \\
&=&2^{\left( s_{1}-\alpha _{1}-n/s\right) N}\left\Vert \eta \right\Vert _{%
\dot{K}_{s}^{\alpha _{1},p}}.
\end{eqnarray*}%
The same arguments give%
\begin{equation*}
\left\Vert f_{N}\right\Vert _{\dot{K}_{s}^{\alpha _{1},p}F_{\theta
}^{s_{1}}}=2^{\left( s_{1}-\alpha _{1}-n/s\right) N}\left\Vert \eta
\right\Vert _{\dot{K}_{s}^{\alpha _{1},p}},\text{ \ \ \ }\left\Vert
f_{N}\right\Vert _{\dot{K}_{q}^{\alpha _{2},r}F_{\theta }^{s_{2}}}=2^{\left(
s_{2}-\alpha _{2}-n/q\right) N}\left\Vert \eta \right\Vert _{\dot{K}%
_{q}^{\alpha _{2},r}},
\end{equation*}%
and%
\begin{equation*}
\left\Vert f_{N}\right\Vert _{\dot{K}_{q}^{\alpha
_{2},p}B_{p}^{s_{2}}}=2^{\left( s_{2}-\alpha _{2}-n/q\right) N}\left\Vert
\eta \right\Vert _{\dot{K}_{q}^{\alpha _{2},p}}.
\end{equation*}%
If the embeddings \eqref{Jawerth}$\ $and \eqref{Frank}$\ $hold then for any $%
N\in \mathbb{N}$%
\begin{equation*}
2^{\left( s_{1}-s_{2}-\alpha _{1}+\alpha _{2}-n/s+n/q\right) N}\leq c.
\end{equation*}%
Thus, we conclude that $s_{1}-\tfrac{n}{s}-\alpha _{1}\leq s_{2}-\tfrac{n}{q}%
-\alpha _{2}$ must necessarily hold by letting $N\rightarrow +\infty $.\vskip%
5pt

Let now $\omega \in \mathcal{S}\left( \mathbb{R}^{n}\right) $ be a function
such that \textrm{supp}$\mathcal{F}\omega \subset \left\{ \xi \in \mathbb{R}%
^{n}:|\xi |<1\right\} $. For $x\in \mathbb{R}^{n}$ and $N\in \mathbb{Z}%
\backslash \mathbb{N}$ we put $f_{N}\left( x\right) =\omega \left(
2^{N}x\right) $. We have $\omega \in \dot{K}_{s}^{\alpha _{1},p}\cap \dot{K}%
_{q}^{\alpha _{2},r}\cap \dot{K}_{q}^{\alpha _{2},p}$. It easy to see that%
\begin{equation*}
\mathcal{F}^{-1}\phi _{j}\ast f_{N}=\left\{ 
\begin{array}{cc}
f_{N}, & j=0 \\ 
0, & \text{otherwise}.%
\end{array}%
\right.
\end{equation*}%
Hence%
\begin{equation*}
\left\Vert f_{N}\right\Vert _{\dot{K}_{s}^{\alpha _{1},p}B_{\beta
}^{s_{1}}}=\left\Vert f_{N}\mid \right\Vert _{\dot{K}_{s}^{\alpha
_{1},p}}=2^{-\left( \alpha _{1}+n/s\right) N}\left\Vert \omega \right\Vert _{%
\dot{K}_{s}^{\alpha _{1},p}}.
\end{equation*}%
The same arguments give%
\begin{eqnarray*}
\left\Vert f_{N}\right\Vert _{\dot{K}_{s}^{\alpha _{1},p}F_{\theta
}^{s_{1}}} &=&2^{-\left( \alpha _{1}+n/s\right) N}\left\Vert \omega
\right\Vert _{\dot{K}_{s}^{\alpha _{1},p}} \\
\left\Vert f_{N}\right\Vert _{\dot{K}_{q}^{\alpha _{2},r}F_{\theta
}^{s_{2}}} &=&2^{-\left( \alpha _{2}+n/q\right) N}\left\Vert \omega
\right\Vert _{\dot{K}_{q}^{\alpha _{2},r}}
\end{eqnarray*}%
and%
\begin{equation*}
\left\Vert f_{N}\right\Vert _{\dot{K}_{q}^{\alpha
_{2},p}B_{p}^{s_{2}}}=2^{-\left( \alpha _{2}+n/q\right) N}\left\Vert \omega
\right\Vert _{\dot{K}_{q}^{\alpha _{2},p}}.
\end{equation*}%
If the embeddings \eqref{Jawerth}$\ $and \eqref{Frank}$\ $hold then for any $%
N\in \mathbb{Z}\backslash \mathbb{N}$ 
\begin{equation*}
2^{-\left( \alpha _{1}-\alpha _{2}+n/s-n/q\right) N}\leq c.
\end{equation*}%
Thus, we conclude that $\alpha _{2}+\tfrac{n}{q}\geq \alpha _{1}+\tfrac{n}{s}
$ must necessarily hold by letting $N\rightarrow -\infty $.\vskip5pt

\section{Applications}

In this section, we give a simple application of Theorems \ \ref%
{embeddings6.1} and \ref{embeddings6.2}. Let $w$ denote a positive, locally
integrable function and $0<p<\infty $. Then the weighted Lebesgue space $%
L^{p}(\mathbb{R}^{n},w)$ contains all measurable functions such that 
\begin{equation*}
\big\|f\big\|_{L^{p}(\mathbb{R}^{n},w)}=\Big(\int_{\mathbb{R}^{n}}\left\vert
f(x)\right\vert ^{p}w(x)dx\Big)^{1/p}<\infty .
\end{equation*}%
For $\varrho \in \lbrack 1,\infty )$ we denote by $\mathcal{A}_{\varrho }$
the Muckenhoupt class of weights, and $\mathcal{A}_{\infty }=\cup _{\varrho
\geq 1}\mathcal{A}_{\varrho }$. We refer to \cite{GR85} for the general
properties of these classes. Let $w\in \mathcal{A}_{\infty }$, $s\in \mathbb{%
R}$, $0<\beta \leq \infty $ and $0<p<\infty $. We define weighted
Triebel-Lizorkin spaces $F_{p,q}^{s}(w)$ to be the set of all distributions $%
f\in \mathcal{S}^{\prime }(\mathbb{R}^{n})$ such that%
\begin{equation*}
\big\|f\big\|_{F_{p,\beta }^{s}(w)}=\Big\|\Big(\sum\limits_{j=0}^{\infty
}2^{js\beta }\left\vert \mathcal{F}^{-1}\varphi _{j}\ast f\right\vert
^{\beta }\Big)^{1/\beta }\Big\|_{L^{p}(\mathbb{R}^{n},w)}
\end{equation*}%
is finite. In the limiting case $q=\infty $ the usual modification is
required. Also we define weighted Besov spaces $B_{p,q}^{s}(w)$ to be the
set of all distributions $f\in \mathcal{S}^{\prime }(\mathbb{R}^{n})$ such
that%
\begin{equation*}
\big\|f\big\|_{B_{p,\beta }^{s}(w)}=\Big(\sum\limits_{j=0}^{\infty
}2^{js\beta }\big\|\mathcal{F}^{-1}\varphi _{j}\ast f\big\|_{L^{p}(\mathbb{R}%
^{n},w)}^{\beta }\Big)^{1/\beta }
\end{equation*}%
is finite. In the limiting case $q=\infty $ the usual modification is
required. The spaces $F_{p,\beta }^{s}(w)$ and $B_{p,\beta }^{s}(w)$ are
independent of the particular choice of the smooth dyadic resolution of
unity $\{\varphi _{j}\}_{j\in \mathbb{N}_{0}}$ appearing in their
definitions. They are quasi-Banach spaces, Banach spaces for $p,q\geq 1$,
moreover for $w\equiv 1\in \mathcal{A}_{\infty }$ we obtain the usual
(unweighted) Besov and Triebel-Lizorkin spaces. Let $w_{\gamma }$ be a power
weight, i.e., $w_{\gamma }(x)=|x|^{\gamma }$ with $\gamma >-n$. Then in view
of the fact that $L^{p}=\dot{K}_{p}^{0,p}$, we have%
\begin{equation*}
\big\|f\big\|_{A_{p,\beta }^{s}(w_{\gamma })}\approx \big\|f\big\|_{\dot{K}%
_{p}^{\frac{\gamma }{p},p}A_{\beta }^{s}}\text{.}
\end{equation*}

Applying Theorems \ref{embeddings6.1} and \ref{embeddings6.2} in some
particular cases yields the following embeddings.

\begin{corollary}
\textit{Let }$s_{1},s_{2}\in \mathbb{R}$, $0<q<s<\infty $, $0<\beta \leq
\infty $ and $w_{\gamma _{1}}(x)=|x|^{\gamma _{1}}$, $w_{\gamma
_{2}}(x)=|x|^{\gamma _{2}}$, with $\gamma _{1}>-n\ $\textit{and }$\gamma
_{2}>-n$. \textit{We suppose that } 
\begin{equation*}
s_{1}-\frac{n+\gamma _{1}}{s}=s_{2}-\frac{n+\gamma _{2}}{q}
\end{equation*}%
and%
\begin{equation*}
\gamma _{2}/q\geq \gamma _{1}/s.
\end{equation*}%
Then%
\begin{equation*}
F_{q,\beta }^{s_{2}}(w_{\gamma _{2}})\hookrightarrow
B_{s,q}^{s_{1}}(w_{\gamma _{1}})\quad \text{and}\quad
B_{q,s}^{s_{2}}(w_{\gamma _{2}})\hookrightarrow F_{s,\beta
}^{s_{1}}(w_{\gamma _{1}}).
\end{equation*}
\end{corollary}




\begin{thebibliography}{99}
\bibitem{BS88} C. Bennett, R. Sharpley, Interpolation of Operators. Academic
Press, San Diego, 1988.

\bibitem{Drihem13} D. Drihem, Embeddings properties on Herz-type Besov and
Triebel-Lizorkin spaces{\normalsize . \emph{Math. Ineq and Appl.} \textbf{16}%
(2) (2013), 439-460.}

\bibitem{Drihem16} D. Drihem, Sobolev embeddings for Herz-type
Triebel-Lizorkin spaces. arXiv:1502.06417.

\bibitem{Fr86} J. Franke, On the spaces $F_{p,q}^{s}$ of Triebel-Lizorkin
type: pointwise multipliers and spaceson domains. \textit{Math Nachr.} 
\textbf{125} (1986), 29-68.

\bibitem{FJ90} {\normalsize M. Frazier, B. Jawerth, A discrete transform and
decomposition of distribution spaces. \emph{J. Funct. Anal.} \textbf{93}(1)
(1990), 34-170.}

\bibitem{GR85} J. Garc\.{\i}a-Cuerva, J.L. Rubio de Francia, \textit{%
Weighted norm inequalities and related topics}, In: North-Holland
Mathematics Studies, Vol. 116, North-Holland, Amsterdam, 1985.

\bibitem{Ja77} B. Jawerth, Some observations on Besov and Lizorkin-Triebel
spaces. \textit{Math. Scand}. \textbf{40} (1977), 94-104.

\bibitem{K10} {\normalsize H. Kempka, Atomic, molecular and wavelet
decomposition of 2-microlocal Besov and Triebel-Lizorkin spaces with
variable integrability. \textit{Funct. Approx. Comment. Math}. \textbf{43}%
(2) (2010), 171-208.}

\bibitem{HerYang99} {\normalsize E. Hernandez, D. Yang, Interpolation of
Herz-type Hardy spaces and applications. \emph{Math. Nachr.} \textbf{42}
(1998), 564-581.}

\bibitem{LiYang96} {\normalsize X. Li, D. Yang, Boundedness of some
sublinear operators on Herz spaces. \emph{Illinois. J. Math.} \textbf{40 }%
(1996), 484-501.}

\bibitem{LuYang1.95} {\normalsize S. Lu, D. Yang, The local versions of $%
H^{p}(\mathbb{R}^{n})$\ spaces at the origin. \emph{Studia. Math.} \textbf{%
116} (1995), 103-131.}

\bibitem{LuYang2.95} {\normalsize S. Lu, D. Yang, The decomposition of the
weighted Herz spaces on $\mathbb{R}^{n}$\ and its applications. \emph{Sci.
in. China (Ser.A).} \textbf{38} (1995), 147-158.}

\bibitem{LuYang97} {\normalsize S. Lu, D. Yang, Herz-type Sobolev and Bessel
potential spaces and their applications. \emph{Sci in China (Ser.A).} 
\textbf{40} (1997), 113-129.}

\bibitem{TD00} {\normalsize L. Tang, D. Yang, Boundedness of vector-valued
operators on weighted Herz spaces. \emph{Approx. Th. \& its Appl.} \textbf{16%
} (2000), 58-70.}

\bibitem{Triebel83} {\normalsize H. Triebel, \emph{Theory of function spaces}%
. Basel: Birkh\"{a}user, 1983.}

\bibitem{Triebel92} {\normalsize H. Triebel, \emph{Theory of function spaces
II}. Basel: Birkh\"{a}user, 1992.}

\bibitem{Triebel06} {\normalsize H. Triebel, \emph{Theory of function spaces
III}. Basel: Birkh\"{a}user, 2006.}

\bibitem{T11} Y. Tsutsui, The Navier-Stokes equations and weak Herz spaces. 
\textit{Advances in Differential Equations}. \textbf{16} (2011), 1049-1085.

\bibitem{Vybiral08} {\normalsize J. Vyb\'{\i}ral, }A new proof of the
Jawerth-Franke embedding. \textit{Rev. Mat. Complut.} \textbf{21}(1) (2008),
75-82.

\bibitem{XuYang03} J. Xu, D. Yang, Applications of Herz-type
Triebel-Lizorkin spaces. \textit{Acta. Math. Sci (Ser. B)}. \textbf{23}
(2003), 328-338.

\bibitem{XuYang05} {\normalsize J. Xu, D. Yang, Herz-type Triebel-Lizorkin
spaces, I.\emph{\ Acta. Math. Sci (English Ed.)}. \textbf{21}(3) (2005),
643-654.}

\bibitem{Xu03} {\normalsize J. Xu, Some properties on Herz-type Besov spaces
(in chinese). \emph{J. Hunan. Univ (Natural Sci).} \textbf{30}(5) (2003),
75-78.}

\bibitem{Xu04} {\normalsize J. Xu, Pointwise multipliers of Herz-type Besov
spaces and their applications.\emph{\ Appl. Math. }\textbf{17}(1) (2004),
115-121.}

\bibitem{Xu05} {\normalsize J. Xu. Equivalent norms of Herz-type Besov and
Triebel-Lizorkin spaces. \emph{J. Funct. Spaces. Appl.} \textbf{3} (2005),
17-31.}

\bibitem{Xu09} {\normalsize J. Xu, Decompositions of non-homogeneous
Herz-type Besov and Triebel-Lizorkin spaces. \emph{Sci. China. Math. }}%
\textbf{57}(2) (2014), 315-331.{\normalsize \ }
\end{thebibliography}
\end{document}